\documentclass[11pt]{article}

\input{amssym}
\begin{document}
\title{{\bf Group classification for the nonlinear \\heat conductivity
equation}}
\date{}
\author{Ali Mahdipour--Shirayeh \\{\small{\it School of Mathematics, Iran
University of Science and Technology,}}\\
{\small{\it Narmak, Tehran 16846-13114, Iran.}~~E-mail:
\textsf{mahdipour@iust.ac.ir }}}
\maketitle
\begin{abstract} Symmetry properties of the nonlinear heat conductivity equations
of the general form $u_t=\left[E(x,u)u_x\right]_x + H(x,u)$ is
studied. The point symmetry analysis of these equations is
considered as well as an equivalence classification which admits
an extension by one dimension of the principal Lie algebra of the
equation. The invariant solutions of equivalence transformations
and classification of the nonlinear heat conductivity equations
among with additional operators are also given.
\end{abstract}
\section{Introduction}
Symmetry properties of mathematical models of heat conductivity
and diffusion processes \cite{Ib1} are traditionally formulated in
terms of nonlinear differential equations which often envisage us
with difficulties in studying. To solve this problem, symmetry
methods play a key role for finding their exact solutions, similar
solutions and invariants \cite{BK,Do,Pa,PS,VP}.

In this study, we generalize the study of a class of the nonlinear
heat conductivity equations (HCEs) which has been recently studied
in some special cases \cite{Ma,PS,VP}. We deal with the class of
nonlinear heat conductivity equations of the general form
\begin{eqnarray}
{\rm HCE}:\quad u_t=\left[E(x,u)u_x\right]_x + H(x,u),
\label{eq:1}
\end{eqnarray}
in which we assumed that $E ,H$ are sufficiently smooth functions,
$E_x, E_u, H_x, H_u \neq 0$, $u$ is treated as the dimensionless
temperature, $t$ and $x$ are the dimensionless time and space
variables and $E$ is the thermal conductivity.

The Lie point symmetry in linear and nonlinear special cases of
our problem was investigated. For instance the problem for the
case in which $E$ just depends to $u$ and $E_u=0$ corresponding to
the linear case, the case where the class of nonlinear
one-dimensional diffusion equations when $E=E(u), H=0$, the class
of diffusion–-reaction equations when $E=E(u), H=\verb"const"$,
the case which the thermal conductivity is a power function of the
temperature and additional equivalence transformations,
conditional equivalence groups and nonclassical symmetries have
all investigated and listed in Table 1 of \cite{VP}. The point
symmetry group of nonlinear fin equations of class (\ref{eq:1})
were considered in a number of papers. Moreover in \cite{Ma},
preliminary group classification of nonlinear fin equation was
studied in the general form.

Class of Eqs.~(\ref{eq:1}) generalizes a great number of the known
nonlinear second order equations describing various processes in
biology, ecology, physics and chemistry (see \cite{BK} and also
\cite{CSR} with references therein). Nonlinear heat equations in
one or higher dimensions are also studied in literature by using
both symmetry as well as other methods \cite{DV,EQZ} (an account
of some interesting cases is given by Polyanin \cite{Po}). But for
the first time, we generalize the equations described above to the
nonlinear heat conductivity equations in the form (\ref{eq:1}) to
investigate their symmetry properties. In order to determine more
symmetry of HCEs, after finding the point symmetry group, we use
preliminary group classification to find different cases of
one-dimensional extension of the symmetry algebra.

The more general class of HCEs is the nonlinear heat conductivity
equations of the form
\begin{eqnarray}
u_t = F(t, x, u,u_x)u_{xx} + G(t, x, u, u_x ), \label{eq:1-2}
\end{eqnarray}
which admits non-trivial symmetry group. The group classification
of (\ref{eq:1-2}) is presented in some references \cite{BLZ,LS}.
However, since the equivalence group of (\ref{eq:1-2}) is
essentially wider than those for particular cases, the results of
\cite{BLZ,LS} cannot be directly used for symmetry classification
of particular ones. Nevertheless, these results are useful for
finding additional equivalence transformations in the class of our
problem. Therefore in contrast to the above works, in the last two
sections of this paper, we study group classification of
Eq.~(\ref{eq:1}) under equivalence transformations in the general
case. Furthermore, a number of nonlinear invariant models which
have nontrivial invariance algebras are obtain.

From \cite{ZL} we know that if the partial differential equation
possesses non-trivial symmetry, then it is invariant under some
finite-dimensional Lie algebra of differential operators which is
completely determined by its structural constants. In the event
that the maximal algebra of invariance is infinite--dimensional,
then it contains, as a rule, some finite-dimensional Lie algebra.
Also, if there are local non-singular changes of variables which
transform a given differential equation into another, then the
finite-dimensional Lie algebra of invariance of these equations
are isomorphic, and in the group-theoretic analysis of
differential equations such equations are considered to be
equivalent. To realize the group classification, we use of the
proposed approach consists in the implementation of an algorithm
explained and performed in references \cite{BLZ,Ib3,Ov,PNI}. For
this goal, our method is similar to the way of \cite{LH} for the
nonlinear wave equation $u_{tt} = f(x,u)u_{xx} + g(x,u)$.

In the next section, we concern with the problem of finding point
symmetry group of Eq.~(\ref{eq:1}). In reminded sections of the
paper we find some further symmetry properties of Eq.~(\ref{eq:1})
by use of equivalence transformations and an extension by one
dimension of the principal Lie algebra of the equation.
\section{Lie point symmetries}
In this section, our study is based on the method of \cite{Ol} for
Lie infinitesimal method. we are concerning with group
classification of HCEs by the point transformations group.

An equation of class (\ref{eq:1}) is a relation among with the
variables of 2--jet space $J^2({\bf R}^2,{\bf R})$ with (local)
coordinate
\begin{eqnarray}
(t,x,u,u_t,u_x,u_{tt},u_{tx},u_{xx}),\label{var}
\end{eqnarray}
where this coordinate involving independent variables $t,x$ and
dependent variable $u$ and derivatives of $u$ in respect to $t$
and $x$ up to order 2 (each index will indicate the derivation
with respect to it, unless we specially state otherwise). Let
${\mathcal M}$ be the total space of independent and dependent
variables resp. $t,x$ and $u$. The solution space of
Eq.~(\ref{eq:1}), (if it exists) is a subvariety
$S_{\Delta}\subset J^2({\bf R}^2,{\bf R})$ of the second order jet
bundle of 2-dimensional sub-manifolds of ${\mathcal M}$. Point
symmetry group on ${\mathcal M}$ is introduced by transformations
in the form of
\begin{eqnarray}
\tilde{t}=\Theta(t,x,u),\hspace{1cm}\tilde{x}=\Xi(t,x,u),\hspace{1cm}
\tilde{u}=\Omega(t,x,u),\label{eq:2}
\end{eqnarray}
for arbitrary smooth functions $\varphi,\chi,\psi$. Also assume
that the general form of infinitesimal generators is
\begin{eqnarray}
Y:=\xi(t,x,u)\,\displaystyle{\frac{\partial }{\partial
t}}+\tau(t,x,u)\,\displaystyle{\frac{\partial }{\partial x}}+
\varphi(t,x,u)\,\displaystyle{\frac{\partial }{\partial u}},
\label{eq:2-1}
\end{eqnarray}
when coefficients are arbitrary smooth functions. These
infinitesimals signify the Lie algebra ${\mathcal L}$ of the point
symmetry group $G$ of Eq.~(\ref{eq:1}). The second order
prolongation of $X$ \cite{Ol,Ov} as a vector field on $J^2({\bf
R}^2,{\bf R})$ is as follows
\begin{eqnarray}
Y^{(2)}:=v + \varphi^t\,\displaystyle{\frac{\partial }{\partial
u_t}}+\varphi^x\,\displaystyle{\frac{\partial }{\partial u_x}}+
\varphi^{tt}\,\displaystyle{\frac{\partial }{\partial u_{tt}}}+
\varphi^{tx}\,\displaystyle{\frac{\partial }{\partial u_{tx}}}+
\varphi^{xx}\,\displaystyle{\frac{\partial }{\partial u_{xx}}},
\end{eqnarray}
where $\eta^t,\eta^x$ and $\eta^{tt},\eta^{tx},\eta^{xx}$ are
arbitrary smooth functions depend to variables $t,x,u,u_t,u_x$ and
(\ref{var}) resp. These coefficients are introduced as following
\begin{eqnarray}
\eta^J &=& {\mathcal D}_J(Q) + \xi\,u_{J,t}+ \tau\,u_{J,x}
\end{eqnarray}
where ${\mathcal D}$ is total derivative, $J$ is a multi-index
with length $1\leq |J|\leq2$ of variables $t,x$ and $Q= u -
\xi\,u_t - \tau\,u_x$ is characteristic of $v$ \cite{Ol}.
According to \cite{Ol}, $v$ is a point infinitesimal generator of
Eq.~(\ref{eq:1}) if and only if
$Y^{(2)}[\mbox{Eq.~(\ref{eq:1})}]=0$ when Eq.~(\ref{eq:1}) is
hold. By applying $Y^{(2)}$ on the equation we have the following
equation
\begin{eqnarray}
&&\tau[E_x\,u_{xx}+E_u\,u_x\,u_{xx}+E_{xx}\,u_x^2+E_{uu}\,u_x^3+E_{xx}\,u_x+E_{ux}\,u_x^2+H_u\,u_x+H_x]
+\varphi[E_u\,u_{xx}+E_{uu}\,u_x^2 \nonumber\\
&& +E_{ux}\,u_x+H_u]
-\varphi^t+\varphi^x[2E_u\,u_x+E_x]+E\,\varphi^{xx}
=0,\hspace{0.5cm}\mbox{whenever}\hspace{0.5cm}
\mbox{Eq.~(\ref{eq:1}) is satisfied}.
\end{eqnarray}

In the extended form of the latter equation when we consider
$u_t=\left[E(x,u)u_x\right]_x + H(x,u)$, functions $\xi,\tau$ and
$\varphi$ only depend to $t, x, u$ rather than other variables
$u_t, u_x, u_{tt}, u_{tx}$ and $u_{xx}$, hence the equation will
be satisfied if and only if the individual coefficients of the
powers of $u_t,u_x$ and their multiplications vanish. This tends
to the following over-determined system of determining equations
\begin{eqnarray}
&&E\,\xi_x=0, \hspace{1.38cm}   E\,\xi_u=0, \hspace{1.38cm} E_u\,\tau+2E\,\tau_u=0,  \nonumber\\
&&E\,\tau_{uu}+E_u\,\tau_u-E_{uu}\,\tau=0, \hspace{1.5cm}
E_x\,\tau+E_u\,\varphi-2E\,\tau_x+E\,\xi_t=0,\nonumber\\
&&H_x(\tau+\varphi)-\varphi_t+E_x\,\varphi_x+E\,\varphi_{xx}+H(\xi_t-\varphi_u)=0,\nonumber\\
&&E_u\,\xi_t+E(\varphi_{uu}-2\tau_{ux})-2E_u\,\tau_x+E_u\,\varphi_u+E_{uu}\,\varphi+2E_{ux}\,\tau=0, \nonumber\\
&&(E_{xx}+2H_u)\tau+E_{ux}\,\varphi+\tau_t+2E_u\,\varphi_x+E_x(\xi_t-\tau_x)+E(2\,\varphi_{ux}-\tau_{xx}).
\label{eq:4}
\end{eqnarray}

The general solution to differential equations (\ref{eq:4}) for
$\xi,\tau$ and $\varphi$ is
\begin{eqnarray}
\xi(t,x,u)=c, \hspace{1cm} \tau(t,x,u)=0, \hspace{1cm}
\varphi(t,x,u)=0,
\end{eqnarray}
for arbitrary constant $c$.
\paragraph{Theorem 1.}{\em A complete set of all infinitesimal generators of the HCE
(\ref{eq:1}) up to point transformations admits the structure of
one-dimensional Lie algebra ${\mathcal L}=\langle
\frac{\partial}{\partial t}\rangle$.}
\\

It is well--known that the existence of a non--fiber--preserving
symmetry usually indicates that one can significantly simplify the
equation by some kind of hodograph--like transformation
interchanging the independent and dependent variables. Since the
kernel of maximal Lie algebra of the HCE (\ref{eq:1}) is
${\mathcal L}=\langle \frac{\partial}{\partial t}\rangle$, thus we
can not use this advantage for simplifying HCEs.

Furthermore, according to the statements of page 209 of \cite{Ol},
we conclude that

\paragraph{Corollary 2.}{\em A system of HCEs of class (\ref{eq:1}) can not be reduced
into an inhomogeneous form of a linear system.}
\section{Equivalence transformations}
In this section, we follow the method of Ovsiannikov \cite{Ov} for
partial differential equations. His approach is based on the
concept of an equivalence group, which is a Lie transformation
group acting in the extended space of independent variables,
functions and their derivatives, and preserving the class of
partial differential equations under study. It is possible to
modify Lie's algorithm in order to make it applicable for the
computation of this group \cite{BLZ,Ib3,Ov,PNI}. Next we construct
the optimal system of subgroups of the equivalence group.

An {\it equivalence transformation} is a non-degenerate change of
the variables $t, x, u$ taking any equation of the form
(\ref{eq:1}) into an equation of the same form, generally
speaking, with different $E(x,u)$ and $H(x,u)$. The set of all
equivalence transformations forms an equivalence group $G$. We
shall find a continuous subgroup $G_C$ of it making use of the
infinitesimal method.

We investigate for an operator of the group $G_C$ in the general
form
\begin{eqnarray}
Y:=\xi(t,x,u)\,\displaystyle{\frac{\partial }{\partial
t}}+\tau(t,x,u)\,\displaystyle{\frac{\partial }{\partial x}}+
\varphi(t,x,u)\,\displaystyle{\frac{\partial }{\partial u}} +
\chi(t,x,u,E,h)\,\displaystyle{\frac{\partial }{\partial E}} +
\eta(t,x,u,E,h)\,\displaystyle{\frac{\partial }{\partial h}}.
\label{eq:10}
\end{eqnarray}
from the invariance conditions of Eq.~(\ref{eq:1}) written as the
system
\begin{eqnarray}
u_t = \left[E(x,u)u_x\right]_x + H(x,u), \hspace{1cm}  E_t = H_t
=0, \label{eq:11}
\end{eqnarray}
where we assumed that $u, E, H$ are differential variables: $u$ on
the base space $(t, x)$ and $E , H$ on the total space $(t, x,
u)$. Also in Eq.~(\ref{eq:10}) the coefficients are dependent to
$t, x, u$ and the two last ones, in addition, depend to $E, h$.
The invariance conditions of the system (\ref{eq:11}) are
\begin{eqnarray}
\widetilde{Y}\,[u_t - \left[E(x,u)u_x\right]_x - H(x,u)]=0,
\hspace{1cm} \widetilde{Y}\,[E_t]=\widetilde{Y}\,[H_t] = 0,
\label{eq:12}
\end{eqnarray}
where
\begin{eqnarray}
\widetilde{Y}:= Y +  \varphi^t\,\displaystyle{\frac{\partial
}{\partial u_t}}+\varphi^x\,\displaystyle{\frac{\partial
}{\partial u_x}}
%+ \varphi^{tt}\,\displaystyle{\frac{\partial
%}{\partial
%u_{tt}}}+\varphi^{tx}\,\displaystyle{\frac{\partial}{\partial
%u_{tx}}}
+ \varphi^{xx}\,\displaystyle{\frac{\partial }{\partial
u_{xx}}}+ \chi^t\,\displaystyle{\frac{\partial }{\partial E_t}} +
\chi^x\,\displaystyle{\frac{\partial }{\partial E_x}} +
\chi^u\,\displaystyle{\frac{\partial }{\partial E_u}}+
\eta^t\,\displaystyle{\frac{\partial }{\partial
H_t}}.\label{eq:12-1}
\end{eqnarray}
is the prolongation of the operator (\ref{eq:10}). Coefficients
$\eta^J$ for multi--index $J$ (with length $1\leq |J|\leq 2$) have
given in section 2 and by applying the prolongation procedure to
differential variables $E, H$ with independent variables $(t, x,
u)$ we have
\begin{eqnarray}
&& \chi^I = \widetilde{{\mathcal D}}_I(\chi) -
E_t\,\widetilde{{\mathcal D}}_I(\xi)- E_x\,\widetilde{{\mathcal
D}}_I(\tau) - E_u\,\widetilde{{\mathcal D}}_I(\phi)=
\widetilde{{\mathcal D}}_I(\chi) - E_x\,\widetilde{{\mathcal
D}}_I(\tau)- E_u\,\widetilde{{\mathcal D}}_I(\phi),  \nonumber\\
[-2.5mm]
&& \\[-2mm]
&& \eta^t = \widetilde{{\mathcal D}}_t(\eta) -
H_t\,\widetilde{{\mathcal D}}_t(\xi)- H_x\,\widetilde{{\mathcal
D}}_t(\eta) - H_u\,\widetilde{{\mathcal D}}_t(\phi) =
\widetilde{{\mathcal D}}_t(\eta)  - H_x\,\widetilde{{\mathcal
D}}_t(\tau) - H_u\,\widetilde{{\mathcal D}}_t(\phi),  \nonumber
\end{eqnarray}
where $I$ varies on variables $t,x,u$ and
\begin{eqnarray}
&& \widetilde{{\mathcal D}}_I := \displaystyle{\frac{\partial
}{\partial I}} + E_I\,\displaystyle{\frac{\partial }{\partial E}}
+ H_I\,\displaystyle{\frac{\partial}{\partial H}}.\nonumber
\end{eqnarray}
Substituting (\ref{eq:12-1}) in (\ref{eq:12}) we tend to the
following system
\begin{eqnarray}
&& \chi\,u_{xx}+\eta - \varphi^t + [2E_u\,u_x+E_x]\,\varphi^x + E\,\varphi^{xx}+\chi^x\,u_x + \chi^u\,u_x^2=0,  \label{eq:13}\\
&& \chi^t=0,  \hspace{1cm}\eta^t=0. \label{eq:14}
\end{eqnarray}

Replacing relations $\eta^J$ (for multi--index $J$ with length
$1\leq|J|\leq 2$) and $\chi^t, \chi^x$ in
Eqs.~(\ref{eq:13})--(\ref{eq:14}) and then introducing the
relation $u_t= \left[E(x,u)u_x\right]_x + H(x,u)$ to eliminate
$u_t$, we have three relations which are called determining
equations. The two last ones are the determining equations
associated with Eqs.~ (\ref{eq:14}), i.e.,
\begin{eqnarray}
&& \chi_t-E_x\,\tau_t-E_u\,\varphi_t=0, \hspace{1cm}
\eta_t-H_x\,\tau_t-H_u\,\varphi_t=0.  \label{eq:15}
\end{eqnarray}
But these relations must hold for arbitrary variables $E_x, E_u,
H_x, H_u$ of the jet space and this fact results in the following
conditions
\begin{eqnarray}
\tau_t =0,  \hspace{1cm}  \varphi_t=0, \hspace{1cm} \chi_t=0,
\hspace{1cm} \eta_t =0,
\end{eqnarray}
so, we find that
\begin{eqnarray}
\xi=\xi(t,x,u),  \hspace{0.7cm}  \eta=\eta(x,u), \hspace{0.7cm}
\varphi=\varphi(x,u), \hspace{0.7cm} \chi=\chi(x,u,E,H),
\hspace{0.7cm} \eta=\eta(x,u,E,H).
\end{eqnarray}
But adding these conditions to the first determining equation,
knowing that $u_t, u_x, u_{tt}, u_{tx}, u_{xx}$ are considered to
be independent variables, we lead to the following system of
equations
\begin{eqnarray}
&&\chi + E\,[\xi_t - 2\,\tau_x]=0, \hspace{1cm}
E\,\varphi_{uu} + \chi_u =0, \hspace{1cm} \xi_t-2\tau_x+\chi_E=0,\nonumber\\
&& E\,[2\,\varphi_{ux}-\tau_{xx}]+\chi_x =0, \hspace{1cm}\varphi_x=0, \hspace{1cm} \chi_H =0,\\
&& \eta  + H\,[\xi_t-\varphi_u] + E\,\varphi_{xx}=0, \hspace{1cm}
E\,\xi_x =0, \hspace{1cm} E\,\xi_u = 0, \hspace{1cm}
E\,\tau_u=0.\nonumber
\end{eqnarray}
This system follows
\begin{eqnarray}
&& \xi = (2\,c_1+c_2)\,t + c_3, \hspace{1cm}\tau = c_1\,x +
c_4, \hspace{1cm} \varphi=c_5\,u + c_6, \nonumber\\[-2.5mm]
&& \\[-2.5mm]
&& \chi=c_2\,E, \hspace{1.1cm} \eta = (c_5-2\,c_1-c_2)\,H,
\nonumber
\end{eqnarray}
for arbitrary constants $c_1,\cdots, c_6$. Therefore the class of
Eqs.~(\ref{eq:1}) has an infinite continuous group of equivalence
transformations generated by infinitesimal operators
\begin{eqnarray}
&& Y_1=\displaystyle{\frac{\partial }{\partial t}}, \hspace{1cm}
Y_2=\displaystyle{\frac{\partial }{\partial x}}, \hspace{1cm}
Y_3=\displaystyle{\frac{\partial }{\partial u}},  \hspace{1cm}
Y_4=2\,t\,\displaystyle{\frac{\partial }{\partial t}+
\frac{\partial }{\partial x}-2\,H\,\frac{\partial }{\partial
H}},\nonumber\\[-2.5mm]
&& \label{eq:15-1}\\[-2.5mm]
&& Y_5=-t\,\displaystyle{\frac{\partial }{\partial t}+
E\,\frac{\partial }{\partial E}+H\,\frac{\partial }{\partial H}},
\hspace{1cm} Y_6=\displaystyle{u\,\frac{\partial }{\partial u}+
H\,\frac{\partial }{\partial H}}.\nonumber
\end{eqnarray}
Moreover, in the group of equivalence transformations are included
also discrete transformations, i.e., reflections
\begin{eqnarray}
t\longmapsto -t,
%\hspace{1cm} x \longmapsto -x,
\hspace{1cm} u \longmapsto -u, \hspace{1cm} E\longmapsto -E,
\hspace{1cm} H\longmapsto -H.\label{eq:1-1-1}
\end{eqnarray}

The communication relations between these vector fields is given
in Table 1. The Lie algebra ${\mathcal L}:=\langle\,  Y_i:
i=1,\cdots,6\,\rangle$ is solvable since the descending sequence
of derived subalgebras of ${\mathcal L}$: ${\mathcal
L}\supset{\mathcal L}^{(1)}=\langle\, Y_1, Y_2, Y_3 \,
\rangle\supset{\mathcal L}^{(2)}=\{0\}$, terminates with a null
ideal. But for each $v=\sum_i v_i Y_i$ and $w=\sum_j w_j Y_j$ in
${\mathcal L}$, its Killing form:
\begin{eqnarray}
K(v,w)={\rm tr}({\rm ad}(v)\circ{\rm
ad}(w))=5\,a_4\,b_4-2\,(a_4\,b_5+a_5\,b_4)+a_5\,b_5+a_6\,b_6,
\end{eqnarray}
is degenerate. Hence ${\mathcal L}$ is neither semisimple nor
simple.\\
\paragraph{Theorem 3.}{\em Let $G_i$ be the one--parameter group (flow) generated by
$Y_i$, then we have
\begin{eqnarray}
\begin{array}{llll}
G_1: (t,x,u,E,h)\longmapsto (t+s,x,u,E,H), &&& \hspace{-0.7cm}
G_2: (t,x,u,E,h)\longmapsto (t,x+s,u,E,H), \\
G_3: (t,x,u,E,h)\longmapsto (t, x, u+s, E, H),
&&&\hspace{-0.7cm} G_4: (t,x,u,E,h)\longmapsto (t\,e^{2s},x\,e^s,u,E,H\,e^{-2s}),\\
G_5: (t,x,u,E,h)\longmapsto (t\,e^{-s}, x, u, E\,e^s, H\,e^s),
&&&\hspace{-0.7cm} G_6: (t,x,u,E,h)\longmapsto
(t,x,u\,e^s,E,H\,e^s),
\end{array}
\end{eqnarray}
when $s$ is an arbitrary parameter. Moreover, if $u=f(t,x)$ for
functions $E$ and $H$ be a solution of the HCE (\ref{eq:1}), so
are
\begin{eqnarray}
&& u_1=f(t+s,x),\hspace{1cm} u_2=f(t,x+s), \hspace{1cm}
u_3=f(t,x,u)-s
\end{eqnarray}
for the same functions $E$ and $H$, $u_4=f(t\,e^{2s},x\,e^s)$ for
the same functions $E$ and $\bar{H}=H\,e^{-2s}$, $u_5=f(t\,e^{-s},
x)$ for $\bar{E}=E\,e^s$ and $\bar{H}=H\,e^s$ and also
$u_6=e^{-s}\,f(t,x)$ for the same $E$ and $\bar{H}=H\,e^s$.}
\section{Preliminary group classification}
In many applications of group analysis, most of extensions of the
principal Lie algebra admitted by the equation under consideration
are taken from the equivalence algebra ${\mathcal L}_{{\mathcal
E}}$. These extensions are called ${\mathcal E}$--extensions of
the principal Lie algebra. The classification of all nonequivalent
equations (with respect to a given equivalence group $G_{\mathcal
E}$) admitting ${\mathcal E}$--extensions of the principal Lie
algebra is called a {\it preliminary group classification}
\cite{Ib3}. We consider the algebra ${\mathcal L}$ spanned on
operators (\ref{eq:15-1}) and use it for a preliminary group
classification.

It is well-known that the problem of classifying invariant
solutions is equivalent to the problem of classifying subgroups of
the full symmetry group under conjugation in which itself is
equivalent to determining all conjugate subalgebras \cite{Ol,Ov}.
The latter problem, tends to determine a list (that is called an
{\it optimal system}) of conjugacy inequivalent subalgebras with
the property that any other subalgebra is equivalent to a unique
member of the list under some element of the adjoint
representation i.e. $\overline{\mathcal L}_H\,{\rm
Ad}(g)\,{\mathcal L}_H$ for some $g$ of a considered Lie group.
Thus we will deal with the construction of the optimal system of
subalgebras of ${\mathcal L}$.

The adjoint action is given by the Lie series
\begin{eqnarray}
{\rm Ad}(\exp(s\,Y_i))Y_j
=Y_j-s\,[Y_i,Y_j]+\frac{s^2}{2}\,[Y_i,[Y_i,Y_j]]-\cdots,
\end{eqnarray}
where $s$ is a parameter and $i,j=1,\cdots,6$. The adjoint
representations of ${\mathcal L}$ is listed in Tables 2; it
consists the separate adjoint actions of each element of
${\mathcal L}$ on all other elements.
\begin{table}[h]
\centering{\caption{Commutator table}}\label{table:1}
\vspace{-0.3cm}
\begin{eqnarray*}\hspace{-0.75cm}\begin{array}{lllllll} \hline
  [\,,\,]  &\hspace{1cm} Y_1     &\hspace{1cm} Y_2 &\hspace{1cm} Y_3  &\hspace{1cm} Y_4   &\hspace{1cm}Y_5 &\hspace{1cm} Y_6 \\ \hline
  Y_1      &\hspace{1cm} 0       &\hspace{1cm} 0   &\hspace{1cm} 0    &\hspace{1cm}2\,Y_1 &\hspace{1cm}-Y_1&\hspace{1cm} 0   \\
  Y_2      &\hspace{1cm} 0       &\hspace{1cm} 0   &\hspace{1cm} 0    &\hspace{1cm}Y_2    &\hspace{1cm}0   &\hspace{1cm} 0   \\
  Y_3      &\hspace{1cm} 0       &\hspace{1cm} 0   &\hspace{1cm} 0    &\hspace{1cm} 0     &\hspace{1cm}0   &\hspace{1cm}Y_3  \\
  Y_4      &\hspace{1cm} -2\,Y_1 &\hspace{1cm}-Y_2 &\hspace{1cm} 0    &\hspace{1cm} 0     &\hspace{1cm}0   &\hspace{1cm} 0   \\
  Y_5      &\hspace{1cm} Y_1     &\hspace{1cm} 0   &\hspace{1cm} 0    &\hspace{1cm} 0     &\hspace{1cm}0   &\hspace{1cm} 0   \\
  Y_6      &\hspace{1cm} 0       &\hspace{1cm} 0   &\hspace{1cm} -Y_3 &\hspace{1cm} 0     &\hspace{1cm}0   &\hspace{1cm} 0   \\
  \hline
\end{array}\end{eqnarray*}
\end{table}
\begin{table}[h]
\centering{\caption{Adjoint table}}\label{table:2}
\vspace{-0.35cm}\begin{eqnarray*}\hspace{0cm}\begin{array}{lllllll}
\hline
  {\rm Ad} &\hspace{1cm} Y_1       &\hspace{1cm} Y_2    &\hspace{1cm} Y_3    &\hspace{1cm} Y_4       &\hspace{1cm}Y_5       &\hspace{1cm} Y_6     \\ \hline
  Y_1      &\hspace{1cm} Y_1       &\hspace{1cm} Y_2    &\hspace{1cm} Y_3    &\hspace{1cm}Y_4-2s\,Y_1&\hspace{1cm}Y_5+s\,Y_1&\hspace{1cm} Y_6     \\
  Y_2      &\hspace{1cm} Y_1       &\hspace{1cm} Y_2    &\hspace{1cm} Y_3    &\hspace{1cm}Y_4-s\,Y_2 &\hspace{1cm}Y_5       &\hspace{1cm} Y_6      \\
  Y_3      &\hspace{1cm} Y_1       &\hspace{1cm} Y_2    &\hspace{1cm} Y_3    &\hspace{1cm} Y_4       &\hspace{1cm}Y_5       &\hspace{1cm}Y_6-s\,Y_3\\
  Y_4      &\hspace{1cm}e^{2s}\,Y_1&\hspace{1cm}e^s\,Y_2&\hspace{1cm} Y_3    &\hspace{1cm} Y_4       &\hspace{1cm}Y_5       &\hspace{1cm} Y_6      \\
  Y_5      &\hspace{1cm}e^{-s}\,Y_1&\hspace{1cm} Y_2    &\hspace{1cm} Y_3    &\hspace{1cm} Y_4       &\hspace{1cm}Y_5       &\hspace{1cm} Y_6      \\
  Y_6      &\hspace{1cm} Y_1       &\hspace{1cm} Y_2    &\hspace{1cm}e^s\,Y_3&\hspace{1cm} Y_4       &\hspace{1cm}Y_5       &\hspace{1cm} Y_6      \\
  \hline
 \end{array}\end{eqnarray*}
\end{table}
\paragraph{Theorem 4.}{\em An optimal system of one-dimensional Lie subalgebras of the
nonlinear HCE (\ref{eq:1}) is provided by those generated by
\begin{eqnarray}
&&\hspace{-0.7cm} 1)\hspace{0.2cm}A^1=Y_1=\partial_t,  \nonumber  \\
&&\hspace{-0.7cm} 2)\hspace{0.2cm}A^2=Y_2=\partial_x,  \nonumber\\
&&\hspace{-0.7cm} 3)\hspace{0.2cm} A^3=Y_3=\partial_u, \nonumber\\
&&\hspace{-0.7cm} 4)\hspace{0.2cm}A^4=Y_4=2t\,\partial_t+\partial_x-2H\,\partial_H,  \nonumber\\
&&\hspace{-0.7cm} 5)\hspace{0.2cm}A^5=Y_5=-t\,\partial_t+E\,\partial_E+H\,\partial_H,  \nonumber\\
&&\hspace{-0.7cm} 6)\hspace{0.2cm} A^6=Y_6=u\,\partial_u + H\,\partial_H, \nonumber\\
&&\hspace{-0.7cm} 7)\hspace{0.2cm}A^7=\pm Y_1+Y_2=\pm \partial_t+\partial_x,  \nonumber \\
&&\hspace{-0.7cm} 8)\hspace{0.2cm}A^8=\pm Y_1+Y_3=\pm \partial_t+\partial_u, \nonumber\\
&&\hspace{-0.7cm} 9)\hspace{0.2cm} A^9=\pm Y_1+Y_6=\pm \partial_t+u\,\partial_u + H\,\partial_H,\nonumber \\
&&\hspace{-0.7cm} 10)\hspace{0.2cm}A^{10}=Y_2+Y_3=\partial_x+\partial_u, \nonumber \\
&&\hspace{-0.7cm} 11)\hspace{0.2cm}A^{11}=\pm Y_2+Y_4=2t\,\partial_t\pm \,\partial_x + \partial_x -2H\,\partial_H,  \nonumber\\
&&\hspace{-0.7cm} 12)\hspace{0.2cm} A^{12}=\pm
Y_2+Y_5=-t\,\partial_t \pm
\partial_x +E\,\partial_E +H\,\partial_H,  \nonumber \\
&&\hspace{-0.7cm} 13)\hspace{0.2cm}A^{13}=Y_2+Y_6=\partial_x+u\,\partial_u+H\,\partial_H, \nonumber \\
&&\hspace{-0.7cm} 14)\hspace{0.2cm}A^{14}=\pm Y_3+Y_4=2t\,\partial_t+\partial_x \pm \partial_u-2H\,\partial_H, \nonumber\\
&&\hspace{-0.7cm} 15)\hspace{0.2cm} A^{15}=\pm
Y_3+Y_5=-t\,\partial_t \pm
\partial_u+E\,\partial_E + H\,\partial_H,  \label{eq:16}\\
&&\hspace{-0.7cm} 16)\hspace{0.2cm}A^{16}=\alpha_1\,Y_4+Y_5=(2\alpha_1-1)t\,\partial_t+\alpha_1\,\partial_x+E\,\partial_E-(2\alpha_1-1)H\,\partial_H,  \nonumber  \\
&&\hspace{-0.7cm} 17)\hspace{0.2cm}A^{17}=\alpha_2\,Y_4+Y_6=2\,\alpha_2 t\,\partial_t+\alpha_2\,\partial_x+u\,\partial_u-(2\alpha_2-1)H\,\partial_H,  \nonumber\\
&&\hspace{-0.7cm} 18)\hspace{0.2cm} A^{18}=\beta_1\,Y_5+Y_6=-\beta_1t\,\partial_t+u\,\partial_u+\beta_1E\,\partial_E+(\beta_1+1)H\,\partial_H,\nonumber \\
&&\hspace{-0.7cm} 19)\hspace{0.2cm} A^{19}=\pm Y_1+Y_2+Y_3=\pm \partial_t+\partial_x+\partial_u,\nonumber  \\
&&\hspace{-0.7cm} 20)\hspace{0.2cm} A^{20}=\pm Y_1+Y_2+Y_6=\pm \partial_t+\partial_x+u\,\partial_u+H\,\partial_H,  \nonumber\\
&&\hspace{-0.7cm} 21)\hspace{0.2cm} A^{21}=\pm Y_2 \pm Y_3+Y_4=2t\,\partial_t\pm \,\partial_x+\partial_x \pm \partial_u-2H\,\partial_H,\nonumber \\
&&\hspace{-0.7cm} 22)\hspace{0.2cm} A^{22}=\pm Y_2\pm Y_3+Y_5=-t\,\partial_t \pm \partial_x \pm \partial_u+E\,\partial_E+H\,\partial_H,   \nonumber\\
&&\hspace{-0.7cm} 23)\hspace{0.2cm} A^{23}=\pm
Y_2+\alpha_3\,Y_4+Y_5=(2\alpha_3-1)t\,\partial_t+(\alpha_3 \pm
1)\partial_x+E\,\partial_E-(2\alpha_3-1)H\,\partial_H, \nonumber\\
&&\hspace{-0.7cm}
24)\hspace{0.2cm}A^{24}=Y_2+\alpha_4\,Y_4+Y_6=2\alpha_4t\,\partial_t+(\alpha_4+1)\partial_x+u\,\partial_u
-(2\alpha_4-1)H\,\partial_H, \nonumber\\
&&\hspace{-0.7cm} 25)\hspace{0.2cm}A^{25}=\pm
Y_2+\beta_2\,Y_5+Y_6=
-\beta_2t\,\partial_t\pm\partial_x+u\,\partial_u+\beta_2\,E\partial_E+(\beta_2+1)H\,\partial_H, \nonumber \\
&&\hspace{-0.7cm} 26)\hspace{0.2cm}A^{26}=\pm
Y_3+\alpha_5\,Y_4+Y_5=
(2\alpha_5-1)t\,\partial_t+\alpha_5\,\partial_x \pm
\partial_u + E\,\partial_E-(2\alpha_5-1)H\,\partial_H,  \nonumber\\
&&\hspace{-0.7cm}
27)\hspace{0.2cm}A^{27}=\alpha_6\,Y_4+\beta_3\,Y_5+Y_6=
(2\alpha_6-\beta_3)t\,\partial_t+\alpha_6\,\partial_x
+u\,\partial_u+\beta_3E\,\partial_E-(2\alpha_6-\beta_3-1)H\,\partial_H, \nonumber  \\
&&\hspace{-0.7cm} 28)\hspace{0.2cm}A^{28}=\pm Y_2\pm
Y_3+\alpha_7\,Y_4+Y_5= (2\alpha_7-1)t\,\partial_t+(\alpha_7\pm
1)\,\partial_x \pm \partial_u+E\,\partial_E-(2\alpha_7-1)H\,\partial_H, \nonumber \\
&&\hspace{-0.7cm} 29)\hspace{0.2cm}A^{29}=\pm
Y_2+\alpha_8\,Y_4+\beta_4\,Y_5+Y_6=
(2\alpha_8-\beta_4)t\,\partial_t+(\alpha_8 \pm 1)\,\partial_x
+u\,\partial_u+\beta_4E\,\partial_E-(2\alpha_8-\beta_4-1)H\,\partial_H,\nonumber
\end{eqnarray}
for nonzero constants $\alpha_i,\beta_j$ ($1\leq i\leq8, 1\leq j\leq 4$).}\\

\noindent{\bf Proof.} Let ${\mathcal L}$ is the symmetry algebra
of Eq.~(\ref{eq:1}) with adjoint representation determined in
Table 2 and
\begin{eqnarray}
Y=a_1\,Y_1+a_2\,Y_2+a_3\,Y_3+a_4\,Y_4+a_5\,Y_5+a_6\,Y_6
\end{eqnarray}
is a nonzero vector field of ${\mathcal L}$. We will simplify as
many of the coefficients $a_i$ as possible through proper adjoint
applications on $Y$. We follow our aim in the below easy cases.
\begin{description}
\item[Case 1]
At first, assume that $a_6\neq 0$. Scaling $Y$ if necessary, we
can consider $a_6$ to be 1 and so follow the problem with
\begin{eqnarray}
Y=a_1\,Y_1+a_2\,Y_2+a_3\,Y_3+a_4\,Y_4+a_5\,Y_5+Y_6.
\end{eqnarray}
\item[Case 1a]
According to Table 2 in the case which $a_4\neq 0$, if we act on
$Y$ by ${\rm Ad}(\exp(\frac{a_1}{a_4}\,Y_2))$, the coefficient of
$Y_1$ can be vanished:
\begin{eqnarray}
Y'=a_2\,Y_2+a_3\,Y_3+a_4\,Y_4+a_5\,Y_5+Y_6.
\end{eqnarray}
Then for $a_3\neq 0$ we apply ${\rm Ad}(\exp(a_3\,Y_3))$ on $Y'$
to cancel the coefficient of $Y_3$ (it is automatically hold for
$a_3=0$):
\begin{eqnarray}
Y''=a_2\,Y_2+a_4\,Y_4+a_5\,Y_5+Y_6.
\end{eqnarray}
\item[Case 1a-1]
If in addition $a_2\neq 0$, we can act ${\rm Ad}(\exp(\pm
\ln(\frac{1}{a_2})\,Y_4))$ on $Y''$ to change $a_2$ to $\pm 1$:
\begin{eqnarray}
Y'''=\pm Y_2+a_4\,Y_4+a_5\,Y_5+Y_6.
\end{eqnarray}
Then for the case which $a_5\neq0$, we can not simplify $Y'''$ any
more. This introduce part 29 of the theorem.\\
Also when $a_5=0$ we tend to part 24 of the theorem.
\item[Case 1a-2]
In Case 1a, let we consider $a_4$ of $Y''$ is equal to zero. Then
for
\begin{eqnarray}
Y''''=a_4\,Y_4+a_5\,Y_5+Y_6,
\end{eqnarray}
where $a_5\neq 0$ more simplification is impossible and this
results in part 27, and where $a_5=0$ it results in part 17 of the
theorem.
\item[Case 1b]
Let the coefficient $a_4$ in $Y$ is zero. So we have the following
new form of $Y$
\begin{eqnarray}
\overline{Y}=a_1\,Y_1+a_2\,Y_2+a_3\,Y_3+a_5\,Y_5+Y_6.
\end{eqnarray}
\item[Case 1b-1]
In the case which $a_5\neq 0$, by applying ${\rm
Ad}(\exp(-\frac{a_1}{a_5})\,Y_1))$ on $\overline{Y}$, the
coefficient $a_1$ will be vanished:
\begin{eqnarray}
\overline{\overline{Y}}=a_2\,Y_2+a_3\,Y_3+a_5\,Y_5+Y_6.
\end{eqnarray}
Furthermore, if we act ${\rm Ad}(\exp(\frac{a_3}{a_6})\,Y_3))$ on
the latter form, so the coefficient $a_3$ will be zero. When
$a_2\neq 0$, applying ${\rm Ad}(\exp(\pm \ln\frac{1}{a_2})\,Y_4))$
any one--dimensional subalgebra generated by $Y$ is equivalent to
one generated by $\pm Y_2+a_5\,Y_5 + Y_6$ which introduce part 25
of the theorem for constant $a_5\neq 0$. For the other case which
$a_2=0$, part 18 are given.
\item[Case 1b-2]
In $\overline{Y}$, let $a_5=0$. In this case, by applying ${\rm
Ad}(\exp(\frac{a_3}{a_6})\,Y_6))$  we can make the coefficient of
$Y_3$ equal to zero:
\begin{eqnarray}
\overline{\overline{\overline{Y}}}=a_1\,Y_1+a_2\,Y_2+Y_6.
\end{eqnarray}
\item[Case 1b-2-1]
For $a_2\neq 0$ the action of ${\rm Ad}(\exp(\pm\ln
\frac{1}{a_2})\,Y_4))$ on $\overline{\overline{\overline{Y}}}$
results in the form $\pm\frac{2a_1}{a_2}\,Y_1+Y_2+Y_6$. Then if
$a_1\neq 0$, we can apply ${\rm
Ad}(\exp(\pm\ln\frac{a_2}{2a_1})\,Y_5))$ on it to have the form
$\pm Y_1+Y_2+Y_6$ that is part 20 of the theorem. If we consider
the case in which $a_1=0$, then we tend to part 13.
\item[Case 1b-2-2]
In $\overline{\overline{\overline{Y}}}$ assume that $a_2=0$. In
this case the condition $a_1\neq 0$ and the action of ${\rm
Ad}(\exp(\pm\ln \frac{1}{a_1})\,Y_5))$ show that the simplest form
of $Y$ is similar to part 9, while the condition $a_1=0$ leads to
part 6 of the theorem.\\

%%%%%%
%%%%%%
%
\item[Case 2]
The remaining one--dimensional subalgebras are spanned by vector
fields of the form $Y$ with $a_6=0$.
\item[Case 2a]
If $a_5\neq 0$ then by scaling $Y$, we can assume that $a_5=1$:
\begin{eqnarray}
\hat{Y}=a_1\,Y_1+a_2\,Y_2+a_3\,Y_3+a_4\,Y_4+Y_5.
\end{eqnarray}
\item[Case 2a-1]
Now by the action of ${\rm Ad}(\exp(\frac{a_1}{a_4}\,Y_2))$ on $Y$
where $a_4\neq 0$, we can cancel the coefficient of $Y_1$:
\begin{eqnarray}
\hat{Y}'=a_2\,Y_2+a_3\,Y_3+a_4\,Y_4+Y_5.
\end{eqnarray}
\item[Case 2a-1-1]
Then for $a_3\neq 0$ by applying ${\rm Ad}(\exp(\pm\ln
\frac{1}{a_3}\,Y_6))$ on $\hat{Y}'$ the coefficient of $Y_3$ can
be $\pm 1$. Thus if $a_2\neq 0$ the action of ${\rm
Ad}(\exp(\pm\ln \frac{1}{a_2}\,Y_4))$ recommends part 28 of the
theorem. Moreover, where $a_2=0$ we tend to part 26.
\item[Case 2a-1-2]
In case 2a-1, let $a_3=0$. So for $a_2\neq0$ applying ${\rm
Ad}(\exp(\pm\ln \frac{1}{a_2}\,Y_4))$ on $\bar{Y}'$ one can lead
to part 23 of the theorem, whereas for $a_2=0$, part 16 are
achieved.
\item[Case 2a-2]
Suppose that $a_4$ in $\hat{Y}$ is equal to zero, then for
$a_1\neq0$ by acting ${\rm Ad}(\exp(-\frac{a_1}{a_5}\,Y_5))$ and
also for $a_1=0$ we lead to the following form
\begin{eqnarray}
\hat{Y}''=a_2\,Y_2+a_3\,Y_3+Y_5.
\end{eqnarray}
\item[Case 2a-2-1]
For $a_3\neq 0$ we can act ${\rm Ad}(\exp(\pm\ln
\frac{1}{a_6}\,Y_6))$ on $\hat{Y}''$ to change the coefficient of
$Y_3$ to be equal to either $+1$ or $-1$. Hence for $a_2\neq 0$ by
applying ${\rm Ad}(\exp(\pm\ln \frac{1}{a_2}\,Y_4))$ and for
$a_2=0$ we find part 22 and part 15 resp.
\item[Case 2a-2-2]
Assume that in $\hat{Y}$, $a_3$ is zero, then for $a_2\neq 0$ by
applying ${\rm Ad}(\exp(\pm\ln \frac{1}{a_2}\,Y_4))$ and for
$a_2=0$ we tend to part 12 and part 5 of the theorem resp.\\

\item[Case 2b]
In the situation of Case 2, suppose that $a_5=0$.
\item[Case 2b-1]
If in addition we assume that $a_4\neq 0$, then if necessary we
can let it equal to $1$:
\begin{eqnarray}
\widetilde{Y}=a_1\,Y_1+a_2\,Y_2+a_3\,Y_3+Y_4.
\end{eqnarray}
So for $a_1\neq 0$ by applying ${\rm
Ad}(\exp(\frac{a_1}{a_4}\,Y_2))$ and also for $a_1=0$, we have the
following form for $\widetilde{Y}$
\begin{eqnarray}
\widetilde{Y}'=a_2\,Y_2+a_3\,Y_3+Y_4.
\end{eqnarray}
\item[Case 2b-1-1]
If $a_3\neq 0$, the action of ${\rm Ad}(\pm\ln
\exp(\frac{1}{a_3}\,Y_6))$ shows that we can make the coefficient
of $Y_3$ equal to $\pm 1$. Then for $a_2\neq 0$ by applying ${\rm
Ad}(\exp(\pm\ln \frac{1}{a_2}\,Y_4))$ and for $a_2=0$ resp. part
21 and part 14 of the theorem are given.
\item[Case 2b-1-2]
In $\widetilde{Y}'$, suppose that $a_3=0$. The simplest possible
form of $Y$ is equal to part 11 after taking $a_2\neq 0$ and after
acting ${\rm Ad}(\pm\ln \exp(\frac{1}{a_2}\,Y_4))$. Moreover, when
$a_2=0$ the simplest case is equal to part 4.
\item[Case 2b-2]
With conditions of Case 2b, in addition let $a_4=0$:
\begin{eqnarray}
\dot{Y}=a_1\,Y_1+a_2\,Y_2+a_3\,Y_3.
\end{eqnarray}
\item[Case 2b-2-1]
Consider $a_3 \neq 0$, then by scaling we can make the coefficient
of $\widehat{Y}$ equal to $1$. By assuming $a_2\neq 0$ and
applying ${\rm Ad}(\exp(\pm\ln\frac{1}{a_2}\,Y_4))$ we find the
following form
\begin{eqnarray}
\dot{Y}'=\pm\frac{2a_1}{a_2}\,Y_1+Y_2+Y_3.
\end{eqnarray}
When $a_1\neq 0$ we can reduce it to part 19 after acting ${\rm
Ad}(-\exp(\frac{a_2}{2a_1}\,Y_5))$, whereas for $a_1=0$ we find
part 10. If we change our assumption on $a_2$ and consider it
equal to 0, for $a_1\neq0$ by applying ${\rm Ad}(\pm\ln
\exp(\frac{1}{a_1}\,Y_5))$ and for $a_1=0$ resp. we lead to part 8
and part 3.
\item[Case 2b-2-2]
Finally if in the same conditions with Case 2b-2, we assume that
$a_3$ be zero, by scaling we can make the coefficient of $Y_2$
equal to $1$ if $a_2\neq 0$:
\begin{eqnarray}
\dot{Y}''=a_1\,Y_1+Y_2.
\end{eqnarray}
If $a_1\neq 0$ we can act ${\rm Ad}(\pm\ln \exp(-\ln
\frac{1}{a_1}\,Y_5))$ on $\dot{Y}''$. Then no further
simplification is possible and then $Y$ is reduced to part 7,
while $a_1=0$ suggests part 2. The last case occurs when we change
the condition on $a_1$ to be equal to zero, which recommends part
1 of the theorem.
\end{description}

There is not any more possible case for studying and the proof is
complete (Note that in the group of equivalence transformations
there are included also discrete transformations,
Eq.~(\ref{eq:1-1-1})).\hfill\ $\square$

The coefficients $E, H$ of Eq.~(\ref{eq:1}) depend on the
variables $x,u$. Therefore, we take their optimal system's
projections on the space $(x, u, E , h)$. The nonzero in
$(x,u)-$axis projections of (\ref{eq:16}) are
\begin{eqnarray}
&&\hspace{-3.5cm} 1)\hspace{0.2cm}Z^1=A^2=A^7=\partial_x, \nonumber \\
&&\hspace{-3.5cm} 2)\hspace{0.2cm}Z^2=A^3=A^8=\partial_u, \nonumber \\
&&\hspace{-3.5cm} 3)\hspace{0.2cm}Z^3=A^4=\partial_x-2H\,\partial_H,  \nonumber \\
&&\hspace{-3.5cm} 4)\hspace{0.2cm}Z^4=A^5=E\,\partial_E+H\,\partial_H,  \nonumber \\
&&\hspace{-3.5cm} 5)\hspace{0.2cm}Z^5=A^6=A^9=u\,\partial_u + H\,\partial_H, \nonumber \\
&&\hspace{-3.5cm} 6)\hspace{0.2cm}Z^6=A^{10}=A^{19}=\partial_x+\partial_u,   \nonumber \\
&&\hspace{-3.5cm} 7)\hspace{0.2cm}Z^7=A^{11}=\pm \partial_x + \partial_x-2H\,\partial_H,  \nonumber \\
&&\hspace{-3.5cm} 8)\hspace{0.2cm}Z^8=A^{12}=\pm \partial_x +E\,\partial_E + H\,\partial_H,   \nonumber \\
&&\hspace{-3.5cm} 9)\hspace{0.2cm}Z^9=A^{13}=A^{20}=\partial_x+u\,\partial_u+H\,\partial_H,  \nonumber \\
&&\hspace{-3.5cm} 10)\hspace{0.2cm}Z^{10}=A^{14}=\partial_x \pm \partial_u-2H\,\partial_H, \nonumber \\
&&\hspace{-3.5cm} 11)\hspace{0.2cm}Z^{11}=A^{15}=\pm \partial_u +E\,\partial_E + H\,\partial_H, \nonumber  \\
&&\hspace{-3.5cm} 12)\hspace{0.2cm}Z^{12}=A^{16}=\alpha_1\,\partial_x+E\,\partial_E-(2\alpha_1-1)H\,\partial_H,   \label{eq:16-1} \\
&&\hspace{-3.5cm} 13)\hspace{0.2cm}Z^{13}=A^{17}=\alpha_2\,\partial_x+u\,\partial_u-(2\alpha_2-1)H\,\partial_H,  \nonumber \\
&&\hspace{-3.5cm} 14)\hspace{0.2cm}Z^{14}=A^{18}=u\,\partial_u+\beta_1E\,\partial_E+(\beta_1+1)H\,\partial_H, \nonumber \\
&&\hspace{-3.5cm} 15)\hspace{0.2cm}Z^{15}=A^{21}=\pm \partial_x +\partial_x \pm \partial_u-2H\,\partial_H, \nonumber \\
&&\hspace{-3.5cm} 16)\hspace{0.2cm}Z^{16}=A^{22}=\pm \partial_x \pm \partial_u+E\,\partial_E+H\,\partial_H,  \nonumber  \\
&&\hspace{-3.5cm} 17)\hspace{0.2cm}Z^{17}=A^{23}=(\alpha_3 \pm 1)\partial_x+E\,\partial_E-(2\alpha_3-1)H\,\partial_H, \nonumber \\
&&\hspace{-3.5cm} 18)\hspace{0.2cm}Z^{18}=A^{24}=(\alpha_4+1)\partial_x+u\,\partial_u-(2\alpha_4-1)H\,\partial_H, \nonumber \\
&&\hspace{-3.5cm} 19)\hspace{0.2cm}Z^{19}=A^{25}=\pm \partial_x+u\,\partial_u+\beta_2\,E\partial_E+(\beta_2+1)H\,\partial_H, \nonumber  \\
&&\hspace{-3.5cm} 20)\hspace{0.2cm}Z^{20}=A^{26}=\alpha_5\,\partial_x \pm \partial_u+E\,\partial_E-(2\alpha_5-1)H\,\partial_H,  \nonumber \\
&&\hspace{-3.5cm} 21)\hspace{0.2cm}Z^{21}=A^{27}=\alpha_6\,\partial_x+u\,\partial_u+\beta_3E\,\partial_E-(2\alpha_6-\beta_3-1)H\,\partial_H,  \nonumber  \\
&&\hspace{-3.5cm} 22)\hspace{0.2cm}Z^{22}=A^{28}=(\alpha_7 \pm 1)\,\partial_x \pm \partial_u+E\,\partial_E-(2\alpha_7-1)H\,\partial_H, \nonumber  \\
&&\hspace{-3.5cm} 23)\hspace{0.2cm}Z^{23}=A^{29}=(\alpha_8 \pm
1)\,\partial_x +u\,\partial_u + \beta_4E\,\partial_E -
(2\alpha_8-\beta_4-1)H\,\partial_H,  \nonumber
\end{eqnarray}
From paper 7 of \cite{Ib3} we conclude that
\paragraph{Proposition 5.}{\label{a}}{\em Let ${\mathcal L}_m:=\langle\,Y_i: i = 1,
\cdots, m \,\rangle$ be an $m$--dimensional algebra. Denote by
$A^i\, (i = 1, \cdots, s,\, 0<s\leq m,\, s \in {\Bbb N})$ an
optimal system of one--dimensional subalgebras of ${\mathcal L}_m$
and by $Z^i\, (i = 1, \cdots, t,\, 0<t\leq s,\, t\in {\Bbb N})$
the projections of $A^i$, i.e., $Z^i = {\rm pr}(A^i)$. If
equations
\begin{eqnarray}
F = F(x,u),\hspace{0.75cm} G =G(x,u),
\end{eqnarray}
are invariant with respect to the optimal system $Z^i$ then the
equation
\begin{eqnarray}
u_t = (F(x,u)\,u_x)_x +  G(x,u)\,u,\label{eq:17}
\end{eqnarray}
admits the operators $X^i=$ projection of $A^i$ on $(t, x, u)$.}\\
\paragraph{Proposition 6.}{\label{b}}{\em Let Eq.~(\ref{eq:17}) and the equation
\begin{eqnarray}
u_t = (\overline{F}(x,u)\,u_x)_x +
\overline{G}(x,u),\label{eq:18}
\end{eqnarray}
be constructed according to Proposition \ref{a} via optimal
systems $Z^i$ and $\overline{Z}^i$ resp. If the subalgebras
spanned on the optimal systems $Z^i$ and $\overline{Z}^i$ resp.
are similar in ${\mathcal L}_m$, then Eqs.~(\ref{eq:17}) and
(\ref{eq:18}) are equivalent with respect to the equivalence group
$G_m$ generated by ${\mathcal L}_m$.}\\

Now by applying Propositions \ref{a} and \ref{b} for the optimal
system (\ref{eq:16-1}), we want to find all nonequivalent
equations in the form of Eq.~(\ref{eq:1}) admitting ${\mathcal
E}$--extensions of the principal Lie algebra ${\mathcal
L}_{{\mathcal E}}$, by one dimension, i.e, equations of the form
(\ref{eq:1}) such that they admit, together with the one basic
operator $\frac{\partial}{\partial t}$ of ${\mathcal L}_1$, also a
second operator $X^{(2)}$. In each case which this extension
occurs, we indicate the corresponding coefficients $E , H$ and the
additional operator $X^{(2)}$.

We perform the algorithm passing from operators
$Z^i\,(i=1,\cdots,23)$ to $E, H$ and $X^{(2)}$ via the following
examples.

Let consider the vector field
\begin{eqnarray}
Z^4=E\,\partial_E +H\,\partial_H,
\end{eqnarray}
then the characteristic equation corresponding to $Z^4$ is
\begin{eqnarray}
\frac{dE}{E}=\frac{dH}{H},
\end{eqnarray}
which determines invariants. Invariants can be taken in the
following form
\begin{eqnarray}
I_1=x, \hspace{0.5cm} I_2= u, \hspace{0.5cm} I_3=\frac{H}{E}.
\label{eq:19}
\end{eqnarray}
\begin{table}
\centering{\caption{The result of the
classification}}\label{table:3} \vspace{-0.35cm}{\small
\begin{eqnarray*} \hspace{-0.75cm}\begin{array}{l l l l l l}
\hline
  N  &\hspace{0cm} Z      &\hspace{0.1cm} \mbox{Invariant $\lambda$}  &\hspace{0.1cm}  \mbox{Equation}  &\hspace{0.3cm} \mbox{Additional operator}\,X^{(2)} \\ \hline
  1  &\hspace{0cm} Z^1    &\hspace{0.1cm} u   &\hspace{0.1cm}  u_t=[\Phi\,u_x]_x+\Psi            &\hspace{0.3cm} \frac{\partial}{\partial x},\,\pm\frac{\partial}{\partial t}+\frac{\partial}{\partial x}
  \\[2mm]
  2  &\hspace{0cm} Z^2    &\hspace{0.1cm} x   &\hspace{0.1cm}  u_t=[\Phi\,u_x]_x+\Psi            &\hspace{0.3cm} \frac{\partial}{\partial u},\,\pm\frac{\partial}{\partial t}+\frac{\partial}{\partial u} \\[2mm]
  3  &\hspace{0.2cm} Z^3    &\hspace{0.1cm} u   &\hspace{0.1cm}  u_t=[\Phi\,u_x]_x+e^{-2(x+\Psi)}  &\hspace{0.3cm} 2t\,\frac{\partial}{\partial t}+\frac{\partial}{\partial x} \\[2mm]
  4  &\hspace{0cm} Z^5    &\hspace{0.1cm} x   &\hspace{0.1cm}  u_t=[\Phi\,u_x]_x+e^{\ln u+\Psi}  &\hspace{0.3cm} t\,\frac{\partial}{\partial t} \\[2mm]
  5  &\hspace{0cm} Z^6    &\hspace{0.1cm}u-x  &\hspace{0.1cm}  u_t=[\Phi\,u_x]_x+\Psi            &\hspace{0.3cm} \frac{\partial}{\partial x}+\frac{\partial}{\partial u},\,\pm \frac{\partial}{\partial t}+\frac{\partial}{\partial x}+\frac{\partial}{\partial u} \\[2mm]
  6  &\hspace{0cm} Z^7    &\hspace{0.1cm} u   &\hspace{0.1cm}  u_t=[\Phi\,u_x]_x+e^{-(x+\Psi)}   &\hspace{0.3cm} t\,\frac{\partial}{\partial t},\, t\,\frac{\partial}{\partial t} + \frac{\partial}{\partial x} \\[2mm]
  7  &\hspace{0cm} Z^8    &\hspace{0.1cm} u   &\hspace{0.1cm}  u_t=[e^{\Phi\pm x}\,u_x]_x+e^{\Psi\pm x} &\hspace{0.3cm} \pm t\,\frac{\partial}{\partial t}+\frac{\partial}{\partial x} \\[2mm]
  8  &\hspace{0cm}Z^9   &\hspace{0.1cm} \ln u-x &\hspace{0.1cm}  u_t=[\Phi\,u_x]_x+e^{x+\Psi}      &\hspace{0.3cm} \frac{\partial}{\partial x} + u\,\frac{\partial}{\partial u},\, \pm \frac{\partial}{\partial t} +\frac{\partial}{\partial x} + u\,\frac{\partial}{\partial u} \\[2mm]
  9  &\hspace{0cm}Z^{10}&\hspace{0.1cm} x\pm u&\hspace{0.1cm}  u_t=[\Phi\,u_x]_x+e^{-2(\Psi\pm u)} &\hspace{0.3cm} 2t\,\frac{\partial}{\partial t}+\frac{\partial}{\partial x}\pm \frac{\partial}{\partial u} \\[2mm]
  10 &\hspace{0cm} Z^{11} &\hspace{0.1cm} x   &\hspace{0.1cm}  u_t=[e^{\Phi\pm u}\,u_x]_x+e^{\Psi\pm u} &\hspace{0.3cm} t\,\frac{\partial}{\partial t} \pm \frac{\partial}{\partial u} \\[2mm]
  11 &\hspace{0cm}Z^{12}(\alpha_1\neq1/2) &\hspace{0.1cm} u   &\hspace{0.1cm}  u_t=[e^{\Phi+\frac{1}{\alpha_1}\,x}\,u_x]_x+e^{\Psi+(2-\frac{1}{\alpha_1}\,x)}  &\hspace{0.3cm} (2\alpha_1-1)t\,\frac{\partial}{\partial t}+\alpha_1\,\frac{\partial}{\partial x} \\[2mm]
  12 &\hspace{0cm}Z^{12}(\alpha_1=1/2)    &\hspace{0.1cm} u   &\hspace{0.1cm}  u_t=[e^{\Phi+\frac{1}{\alpha_1}\,x}\,u_x]_x+\Psi                                &\hspace{0.3cm}  \frac{\partial}{\partial x} \\[2mm]
  13 &\hspace{0cm}Z^{13}(\alpha_2\neq1/2) &\hspace{0.1cm} \ln u-\frac{x}{\alpha_2}  &\hspace{0.1cm}  u_t=[\Phi\,u_x]_x+e^{(1-2\alpha_2)(\ln u+\Psi)}        &\hspace{0.3cm} 2\alpha_2\,t\,\frac{\partial}{\partial t}+\alpha_2\,\frac{\partial}{\partial x}+u\,\frac{\partial}{\partial u} \\[2mm]
  14 &\hspace{0cm}Z^{13}(\alpha_2=1/2)    &\hspace{0.1cm} \ln u-2x         &\hspace{0.1cm}  u_t=[\Phi\,u_x]_x+\Psi   &\hspace{0.3cm} 2t\,\frac{\partial}{\partial t}+\frac{\partial}{\partial x}+2u\,\frac{\partial}{\partial u} \\[2mm]
  15 &\hspace{0cm}Z^{14}(\beta_1\neq-1)   &\hspace{0.1cm} x   &\hspace{0.1cm}  u_t=[e^{\!\beta_1(\Phi(u)\!+\!\ln u)}\,u_x]_x\!+\!e^{\!(1+\beta_1)(\Psi\!+\!\ln u)}          &\hspace{0.3cm} \beta_1 t\,\frac{\partial}{\partial t} -u\,\frac{\partial}{\partial u} \\[2mm]
  16 &\hspace{0cm}Z^{14}(\beta_1=-1)      &\hspace{0.1cm} x   &\hspace{0.1cm}  u_t=[u\Phi\,u_x]_x+\Psi     &\hspace{0.3cm} t\,\frac{\partial}{\partial t} - u\,\frac{\partial}{\partial u} \\[2mm]
  17 &\hspace{0cm}Z^{15}(2\partial_x\pm\partial_u-2H\partial_H) &\hspace{0.1cm} \frac{1}{2}\,x\pm u   &\hspace{0.3cm}  u_t=[\Phi(u)\,u_x]_x+e^{-2(\Psi\pm u)}  &\hspace{0.3cm} 2\,t\,\frac{\partial}{\partial t} + 2\,\frac{\partial}{\partial x}\pm\frac{\partial}{\partial u}  \\[2mm]
  18 &\hspace{0cm}Z^{15}(\pm\partial_u-2H\partial_H)      &\hspace{0.1cm} x   &\hspace{0.1cm}  u_t=[\Phi\,u_x]_x+e^{-2(\Psi\pm u)}     &\hspace{0.3cm} 2\,t\,\frac{\partial}{\partial t} \pm\,\frac{\partial}{\partial u} \\[2mm]
  19 &\hspace{0cm}Z^{16}  &\hspace{0.1cm} \pm x\pm u   &\hspace{0.1cm}  u_t=[e^{\Phi\pm u}\,u_x]_x+e^{\Psi\pm u}     &\hspace{0.3cm} -t\,\frac{\partial}{\partial t} \pm \frac{\partial}{\partial x} \pm \frac{\partial}{\partial u} \\[2mm]
  20 &\hspace{0cm}Z^{17}(\alpha_3\neq\pm 1,1/2)    &\hspace{0.1cm} u   &\hspace{0.1cm}  u_t=[e^{\Phi+(\alpha_3\pm1)x}\,u_x]_x+e^{\Psi-\frac{2\alpha_3-1}{\alpha_3\pm1}\,x}     &\hspace{0.3cm} (2\alpha_3-1)t\,\frac{\partial}{\partial t} + (\alpha_3\pm1)\,\frac{\partial}{\partial x} \\[2mm]
  21 &\hspace{0cm}Z^{17}(\alpha_3=1/2)       &\hspace{0.1cm} u         &\hspace{0.1cm}  u_t=[e^{\Phi+(1/2\pm1)x}\,u_x]_x+\Psi     &\hspace{0.3cm}     \frac{\partial}{\partial x} \\[2mm]
  22 &\hspace{0cm}Z^{18}(\alpha_4\neq-1,1/2) &\hspace{0.1cm} \ln u-\frac{x}{\alpha_4+1}   &\hspace{0.3cm}  u_t=[\Phi\,u_x]_x+e^{\Psi-\frac{2\alpha_4-1}{\alpha_4+1}\,x}     &\hspace{0.3cm} 2\alpha_4\,t\,\frac{\partial}{\partial t} + (\alpha_4+1)\,\frac{\partial}{\partial x}+u\,\frac{\partial}{\partial u} \\[2mm]
  23 &\hspace{0cm}Z^{18}(\alpha_4=-1)        &\hspace{0.1cm} x   &\hspace{0.1cm}  u_t=[\Phi\,u_x]_x+e^{3(\Psi+\ln u)}     &\hspace{0.3cm} 2\,t\,\frac{\partial}{\partial t} -u\,\frac{\partial}{\partial u} \\[2mm]
  24 &\hspace{0cm}Z^{18}(\alpha_4=1/2)       &\hspace{0.1cm} \ln u-\frac{2}{3}\,x   &\hspace{0.1cm}  u_t=[\Phi\,u_x]_x+\Psi         &\hspace{0.3cm}    t\,\frac{\partial}{\partial t} +\frac{3}{2}\,\frac{\partial}{\partial x}+u\,\frac{\partial}{\partial u} \\[2mm]
  25 &\hspace{0cm}Z^{19}(\beta_2\neq-1)      &\hspace{0.1cm} \ln u\pm x     &\hspace{0.1cm}  u_t=[e^{\beta_2(\Phi\pm x)}\,u_x]_x+e^{(\beta_2+1)(\Psi\pm x)}     &\hspace{0.3cm} -\beta_2t\,\frac{\partial}{\partial t} \pm \frac{\partial}{\partial x}+u\,\frac{\partial}{\partial u} \\[2mm]
  26 &\hspace{0cm}Z^{19}(\beta_2=-1)         &\hspace{0.1cm} \ln u\pm x     &\hspace{0.1cm}  u_t=[e^{\Phi\pm x}\,u_x]_x+\Psi     &\hspace{0.3cm} t\,\frac{\partial}{\partial t} \pm \frac{\partial}{\partial x}+u\,\frac{\partial}{\partial u} \\[2mm]
  27 &\hspace{0cm}Z^{20}(\alpha_5\neq1/2)    &\hspace{0.1cm}\frac{x}{\alpha_5} \pm u &\hspace{0.1cm}  u_t=[e^{\Phi\pm u}\,u_x]_x+e^{\Psi-(2-\frac{1}{\alpha_5})x}     &\hspace{0.3cm} (2\alpha_5-1)t\,\frac{\partial}{\partial t} + \alpha_5\,\frac{\partial}{\partial x}\pm \frac{\partial}{\partial u} \\[2mm]
  28 &\hspace{0cm}Z^{20}(\alpha_5=1/2)       &\hspace{0.1cm} \frac{x}{2} \pm u      &\hspace{0.1cm}  u_t=[e^{\Phi\pm u}\,u_x]_x+\Psi     &\hspace{0.3cm} \frac{\partial}{\partial x} \pm2\,\frac{\partial}{\partial u} \\[2mm]
  29 &\hspace{0cm}Z^{21}(2\alpha_6-\beta_3\neq1)    &\hspace{0.1cm}\ln u-\frac{x}{\alpha_6}   &\hspace{0.1cm}  u_t=[e^{\beta_3(\Phi+\alpha_6x)}\,u_x]_x+e^{\Psi-\frac{2\alpha_6-\beta_3-1}{\alpha_6}\,x}     &\hspace{0.3cm} (2\alpha_6-\beta_3)t\,\frac{\partial}{\partial t} +\alpha_6\,\frac{\partial}{\partial x}+u\,\frac{\partial}{\partial u} \\[2mm]
  30 &\hspace{0cm}Z^{21}(2\alpha_6-\beta_3 =  1)    &\hspace{0.1cm}\ln u-\frac{x}{\alpha_6}   &\hspace{0.1cm}  u_t=[e^{\beta_3(\Phi+\alpha_6x)}\,u_x]_x+\Psi     &\hspace{0.3cm} t\,\frac{\partial}{\partial t} + \alpha_6\,\frac{\partial}{\partial x}+u\,\frac{\partial}{\partial u} \\[2mm]
  31 &\hspace{0cm}Z^{22}(\alpha_7\neq\pm 1,1/2)  &\hspace{0.1cm} \frac{x}{\alpha_7\pm 1} \pm u  &\hspace{0.1cm}  u_t=[e^{\Phi\pm u}u_x]_x+e^{(1-2\alpha_7)(\Psi\pm u)}     &\hspace{0.3cm} (2\alpha_7\!-\!1)t\,\frac{\partial}{\partial t} + (\alpha_7\!\pm\!1)\,\frac{\partial}{\partial x} \pm \frac{\partial}{\partial u} \\[2mm]
  32 &\hspace{0cm}Z^{22}(\alpha_7      = \pm 1)  &\hspace{0.1cm} x                      &\hspace{0.1cm}  u_t=[e^{\Phi\pm u}\,u_x]_x+e^{(1\pm 2)(\Psi\pm u)}     &\hspace{0.3cm} (\pm 2-1)\,t\,\frac{\partial}{\partial t} \pm \frac{\partial}{\partial u} \\[2mm]
  33 &\hspace{0cm}Z^{22}(\alpha_7      =   1/2)  &\hspace{0.1cm} \frac{x}{1/2\pm 1} \pm u       &\hspace{0.1cm}  u_t=[e^{\Phi\pm u}\,u_x]_x+\Psi            &\hspace{0.3cm} (\frac{1}{2}\pm1)\,\frac{\partial}{\partial x} \pm \frac{\partial}{\partial u} \\[2mm]
  34 &\hspace{0cm}Z^{23}(\alpha_8\neq\pm 1, 2\alpha_8-\beta_4\neq1) &\hspace{0.1cm}\ln u-\frac{x}{\alpha_8\pm 1}   &\hspace{0.1cm}  u_t=[e^{\!\beta_4(\Phi+\!\ln u)}u_x]_x\!+\!e^{\!\Psi-\!\frac{2\alpha_8-\beta_4-1}{\alpha_8\pm1}x}     &\hspace{0.3cm} (2\alpha_8\!-\!\beta_4)t\frac{\partial}{\partial t}\!+\! (\alpha_8\!\pm\!1)\frac{\partial}{\partial x}\!+\!u\frac{\partial}{\partial u} \\[2mm]
  35 &\hspace{0cm}Z^{23}(\alpha_8\neq\pm 1, 2\alpha_8-\beta_4=1)    &\hspace{0.1cm}\ln u-\frac{x}{\alpha_8\pm 1}   &\hspace{0.1cm}  u_t=[e^{\beta_4(\Phi+\ln u)}\,u_x]_x+\Psi     &\hspace{0.3cm} t\,\frac{\partial}{\partial t} + (\alpha_8\pm1)\,\frac{\partial}{\partial x}+u\,\frac{\partial}{\partial u} \\[2mm]
  36 &\hspace{0cm}Z^{23}(\alpha_8 =\pm1, 2\alpha_8-\beta_4\neq1)    &\hspace{0.1cm} x  &\hspace{0.1cm}  u_t=[e^{\!\beta_4(\Phi+\!\ln u)}\,u_x]_x\!+\!e^{\!(2\alpha_8\!-\!\beta_4\!-\!1)(\Psi+\!\ln u)}     &\hspace{0.3cm} (\pm2-\beta_4)\,t\,\frac{\partial}{\partial t} + u\,\frac{\partial}{\partial u} \\[2mm]
  37 &\hspace{0cm}Z^{23}(\alpha_8 =\pm1,    2\alpha_8-\beta_4=1)    &\hspace{0.1cm} x  &\hspace{0.1cm}  u_t=[e^{(\pm2-1)(\Phi+\ln u)}\,u_x]_x+\Psi     &\hspace{0.3cm}  u\,\frac{\partial}{\partial u} \\
  \hline
\end{array}\end{eqnarray*}}
\end{table}
In this case there are no invariant equations because the
necessary condition for existence of invariant solutions (see
\cite{Ov}, Section 19.3) is not satisfied, i.e., invariants
(\ref{eq:19}) cannot be solved with respect to $E$ and $H$ since
each two of them can not be an invariant function with respect to
the third one.

Considering $Z^{23}$ (for $a_8 \pm 1 \neq 0$, $2\alpha_8-\beta_4-1
\neq 0$) we have the below characteristic equation
\begin{eqnarray}
\frac{dx}{\alpha_8 \pm
1}=\frac{du}{u}=\frac{dE}{\beta_4E}=\frac{dH}{(2\alpha_8-\beta_4-1)H},
\end{eqnarray}
This equation suggest the following invariants
\begin{eqnarray}
I_1=\ln u - (\alpha_8 \pm 1)x, \hspace{0.5cm} I_2= \ln
E^{1/\beta_4}-\ln u, \hspace{0.5cm} I_3=\ln H +
\frac{2\alpha_8-\beta_4-1}{\alpha_8 \pm 1}\,x.\label{eq:20}
\end{eqnarray}
From the invariance equations we can write
\begin{eqnarray}
I_2=\Phi(I_1), \hspace{1cm} I_3= \Psi(I_1),
\end{eqnarray}
They result in the forms
\begin{eqnarray}
E=\exp(\beta_4(\ln u + \Phi(\lambda))),\hspace{1cm}
H=\exp\Big(\Psi(\lambda)-\frac{2\alpha_8-\beta_4-1}{\alpha_8 \pm
1}\,x\Big),
\end{eqnarray}
for invariant $\lambda=\ln u - (\alpha_8 \pm 1)x$.

From Proposition \ref{a} applied to the operator $Z^{23}$ (for
$a_8 \pm 1 \neq 0$, $2\alpha_8-\beta_4-1 \neq 0$) we obtain the
additional operator
\begin{eqnarray}
X^{(2)}= (2\alpha_8-\beta_4)t\,\partial_t+(\alpha_8 \pm
1)\,\partial_x +u\,\partial_u.
\end{eqnarray}
\mbox{ }\hspace{3mm}   One can perform the algorithm for other
$Z^i$~s of (\ref{eq:16-1}) similarly. The preliminary group
classification of nonlinear fin equation (\ref{eq:1}) admitting an
extension ${\mathcal L}_2$ of the principal Lie algebra ${\mathcal
L}_1$ is listed in Table 3.
\section{Conclusion}

Symmetry analysis for equations $u_t=\left[E(x,u)u_x\right]_x +
H(x,u)$ rather than previous results on special cases of this
equation equation \cite{PS,VP}, is carried out exhaustively. Also,
equivalence classification is given of the equation admitting an
extension by one of the principal Lie algebra of the equation. The
paper is one of few applications of a new algebraic approach to
the problem of group classification: the method of preliminary
group classification. Derived results are summarized in Table 3.
%
%%%%%%%%%%%%%%%%%%%%%%%%%%%
%

%
\end{document}